\documentclass{article}
\usepackage{latexsym, color}
\usepackage{amsthm}
\usepackage[latin1]{inputenc}
\usepackage{a4wide}
\usepackage{amsthm}
\usepackage{graphicx}
\usepackage{amssymb}
\usepackage{amsmath}
\usepackage{ascmac}
\usepackage{setspace}
\usepackage{float}
\usepackage{algorithm}
\usepackage{algorithmic}
\usepackage{setspace}
\usepackage{comment}
\newtheorem{Theorem}[equation]{Theorem}
\newtheorem{Corollary}[equation]{Corollary}
\newtheorem{Lemma}[equation]{Lemma}

\theoremstyle{definition}
\newtheorem{Definition}[equation]{Definition}

\theoremstyle{remark}
\newtheorem{Remark}[equation]{Remark}
\numberwithin{equation}{section}
\newtheorem{Claim}[equation]{Claim}
\DeclareMathOperator{\Hom}{Hom}

\DeclareMathOperator{\id}{id}

\DeclareMathOperator{\Aut}{Aut}

\numberwithin{equation}{section}
\allowdisplaybreaks
\begin{document}
\title{The Yangian associated with the binary tree}
\author{Mamoru Ueda\thanks{Graduate school of Mathematical Sciences, the University of Tokyo, mueda@ms.u-tokyo.ac.jp}}
\date{}

\maketitle

\begin{abstract}
In this article, we give three candidates for a definition of the Yangian associated with the binary tree $X^*$ and the group $\Aut(X^*)$ consisting of automorphisms of $X^*$ which preserve the tree structure. The first candidate has the Schur-Weyl duality with one kind of an extension of the degenerate affine Hecke algebra. The second and third ones have a homomorphism corresponding to the wreath recursion.
\end{abstract}

\section{Introduction}
For a positive integer $l$, the set $X^*=\{(x_i)_{i=1}^m\mid m\geq0,x_i\in\mathbb{Z}/l\mathbb{Z}\}$ can be regarded as a rooted tree whose root is the empty word $\emptyset$. Let us denote a group consisting of automorphisms of the tree $X^*$ by $\Aut(X^*)$.  Naturally, we can embed the symmetric group $S_l$ of degree $l$ to $\Aut(X^*)$ by restricting it to $\{1,\cdots,l\}\subset X^*$. Weyl constructed the Schur-Weyl duality between the symmetric group and the general linear Lie algebra $\mathfrak{gl}(n)$, which gives the equivalence between the category of finite dimensional $S_l$-modules and the category of finite dimensional representations of $\mathfrak{gl}(n)$ of type $l$. The Schur-Weyl dualities were khown for other algebras. For example, the Schur-Weyl duality between Yangians and degenerate affine Hecke algebras was given by Drinfeld \cite{Dri} and Arakawa \cite{ara}. The degenerate affine Hecke algebra of type $A_{l-1}$ is a deformation of $\mathbb{C}[S_l]\ltimes\mathbb{C}[u_1,\cdots,u_l]$ and comes from the Cherednik-Dunkl operators (see \cite{Ch}). The Yangian was one kind of a quantum group and was introduced by Drinfeld \cite{D1} and \cite{D2} as a tool for solving the quantum Yang-Baxter equation. In this article, we give three candidates of the definition of the Yangian associated with $X^*$ in the case $l=2$. In the case $l=2$, the tree is called the binary tree.

In Section~4, we give the Schur-Weyl duality between $\Aut(X^*)$ and the associative algebra generated by $\{E_{i,j}^\sigma\mid1\leq i,j\leq n,\sigma\in\Aut(X^*)\}$. We denote this associative algebra by $G(n)$. Since the Yangian associated with $\mathfrak{gl}(n)$ contains the universal enveloping algebra of $\mathfrak{gl}(n)$, the Yangian associated with $X^*$ should contain $G(n)$. In Section 6, we construct $\mathbb{H}$ as an extension of the degenerate affine Hecke algebra of type $A_1$ and introduce the associative algebra $YY(n)$, which has the Schur-Weyl duality with $\mathbb{H}$. Considering the Schur-Weyl duality between Yangians and degenerate affine Hecke algebras, $YY(n)$ is one of candidates for the Yangian associated with the binary tree.

For $\sigma\in\Aut(X^*)$ and $r=1,2$, the automorphisms $\sigma|_r\in\Aut(X^*)$ and $\sigma_0\in\Aut(X^*)$ can be defined as the restriction of $\sigma$ to the subsets $X^*|_r=\{(x_i)_{i=1}^m\mid m\geq1,x_i\in\mathbb{Z}/l\mathbb{Z}\}$ and $\{1,\cdots,l\}$. 
The group $\Aut(X^*)$ has an isomorphism $\Aut(X^*)\simeq\mathbb{C}[S_l]\ltimes\Aut(X^*)^{\otimes2}$ given by $\sigma\mapsto(\sigma_0,\sigma|_1\otimes\sigma|_2)$. This isomorphism is called the wreath recursion of $\Aut(X^*)$. We expect that the Yangian associated with $X^*$ has a homomorphism corresponding to the wreath recursion. Unfortunately, by the similar way to the original Yangian, we cannot obtain such homomorphisms for $YY(n)$. In the latter half of this paper, by changing the definition of $YY(n)$ a little, we give two candidates, which have a homomorphism corresponding to the wreath recursion.

Let $X=X^*\sqcup\{\cdot\}$ and $P=\{(a_w)_{w\in X}\mid \sum_{w\in X}a_w<\infty\}$. We can define the restriction of $(a_w)\in P$ naturally and denote them by $(a_w)|_1$ and $(a_w)|_2$. In Section~8, we give the associative algebra $Y(n)$ generated by $\{T^{\mathfrak{p}}_{i,j}\mid1\leq i,j\leq n,\mathfrak{p}\in P\}$ which has a homomorphism
\[\Delta\colon Y(n)\to Y(n)\otimes Y(n)\]
determined by $T^{\mathfrak{p}}_{i,j}\mapsto\sum_{\substack{\mathfrak{r}+\mathfrak{s}=\mathfrak{p}}}T^{\mathfrak{r}|_1}_{i,x}\otimes T^{\mathfrak{s}|_2}_{x,j}$.

Since the Yangian associated with $X^*$ should contain $G(n)$, we extend $Y(n)$ in two ways. The group $\Aut(X^*)$ acts on $Y(n)$ by $\sigma T^{\mathfrak{p}}_{i,j}=T^{\sigma(\mathfrak{p})}_{i,j}$ so that we obtain the wreath product $\widehat{Y}(n)=\Aut(X^*)\ltimes Y(n)$.
\begin{Theorem}
There exists a homomorphism $\Delta^D\colon \widehat{Y}(n)\to\widehat{Y}(n)^{\otimes 4}\rtimes\mathbb{C}[S_2]$.
\end{Theorem}
We can give another homomorphism corresponding to the wreath recursion.
In Section~7, we define a different group structure on $\Aut(X^*)$ and denote this group by $\Aut(X^*)^*$. The group $\Aut(X^*)^*$ acts on $Y(n)$ by $\sigma T^{\mathfrak{p}}_{i,j}=T^{\sigma_0\sigma(\mathfrak{p})}_{i,j}$ so that we obtain the wreath product $\widetilde{Y}(n)=\Aut(X^*)^*\ltimes Y(n)$.
\begin{Theorem}
There exists a homomorphism $\widetilde{\Delta}\colon \widetilde{Y}(n)\to\widehat{Y}(n)^{\otimes 2}$
\end{Theorem}
Actually, $\Aut(X^*)$ and $\Aut(X^*)^\ast$ are contained in $\mathbb{C}[S_2]\ltimes\Aut(X^*)^*$. Also, $\mathbb{C}[S_2]\ltimes\Aut(X^*)^*$ has an action on $Y(n)$ which is compatible with actions of $\Aut(X^*)$ and $\Aut(X^*)^\ast$. We conjecture that $\widetilde{\Delta}$ will be extended to $Y(n)\rtimes(\mathbb{C}[S_2]\ltimes\Aut(X^*)^*)$.

The group $\Aut(X^*)$ is one kind of self-similar groups. We expect that we can define the Yangian associated with other self-similar groups and give results corresponding to the solution of the quantum Yang-Baxter equation similarly to the original Yangian.

\section*{Acknowledgments}
The author expresses his sincere thanks to Shuichi Harako, Masato Shinoda and Taisuke Matsuda for helpful discussions. This research was supported by JSPS Research Fellowship for Young Scientists (PD), Grant Number JP25KJ0038.

\section{Schur-Weyl duality for $\mathfrak{gl}(n)$ and the symmetric group $S_n$}
For a positive integer $n$, let $\mathfrak{gl}(n)$ be the general linear Lie algebra spanned by $\{E_{i,j}\mid 1\leq i,j\leq n\}$ with the commutator relation
\[
[E_{i,j},E_{k,l}]=\delta_{j,k}E_{i,l}-\delta_{i,l}E_{k,j}\text{ for }1\leq i,j,k,l\leq n.
\]
Since we can interpret $\mathfrak{gl}(n)$ as the set of $n\times n$-matrices, $\mathfrak{gl}(n)$ naturally acts on $\mathbb{C}^n=\bigoplus_{i=1}^n\mathbb{C}e_i$, where $\{e_i\}_{i=1}^n$ is the standard basis of $\mathbb{C}^n$.

Let us denote the symmetirc group of degree $n$ by $S_n$ and the permutation of $i$ and $j$ by $\sigma_{i,j}$. 
The following duality between the symmetric group and $\mathfrak{gl}(n)$ is called the Schur-Weyl duality (see Barele-Regev \cite{BR}). 
\begin{Theorem}\label{SW1}
\textup{(1)}\ For positive integers $l$ and $n$, the symmetric group $S_l$ acts on $(\mathbb{C}^n)^{\otimes l}$ by
\begin{equation*}
\sigma(v_1\otimes\cdots\otimes v_l)=v_{\sigma(1)}\otimes\cdots\otimes v_{\sigma(l)}
\end{equation*}
for $v_1,\cdots,v_l\in\mathbb{C}^n$ and $\sigma\in S_l$.

\textup{(2)}\ For positive integers $l$ and $n$, the Lie algebra $\mathfrak{gl}(n)$ acts on $(\mathbb{C}^n)^{\otimes l}$ by
\[E_{i,j}v=\sum_{r=1}^n\limits E^{(r)}_{i,j}v\qquad\text{ for }v\in(\mathbb{C}^n)^{\otimes l},1\leq i,j\leq n,\]
where $E^{(r)}_{i,j}=1^{\otimes r-1}\otimes E_{i,j}\otimes 1^{\otimes l-r}$.

\textup{(3)} Suppose that $v=e_{i_1}\otimes\cdots\otimes e_{i_l}$ for distinct $1\leq i_1,\cdots,i_l\leq n$. Then, $v$ generates $(\mathbb{C}^n)^{\otimes l}$ via the action of $U(\mathfrak{gl}(n))$.

\textup{(4)}\ For positive integers $l,n$ and a right module $M$ over $\mathbb{C}[S_l]$, $\mathfrak{gl}(n)$ acts on $M\otimes_{\mathbb{C}[S_l]}(\mathbb{C}^n)^{\otimes l}$ by
\[
E_{i,j}(m\otimes v)=\sum_{r=1}^l\limits m\otimes E^{(r)}_{i,j}v
\]
for $1\leq i,j\leq n,m\in M,v\in(\mathbb{C}^n)^{\otimes l}.$
\end{Theorem}
We say that a $\mathfrak{gl}(n)$-module $N$ is of level $l$ if $N$ is a direct sum of submodules of $(\mathbb{C}^n)^{\otimes2}$.
\begin{Theorem}\label{SW2}
If $n>l$, the functor
\[
SW\colon M\mapsto M\otimes_{\mathbb{C}[S_l]}(\mathbb{C}^n)^{\otimes l}
\]
gives the equivalence between the category of finite dimensional right modules over $\mathbb{C}[S_l]$ and the category of finite dimensional representations of $\mathfrak{gl}(n)$ of level $l$.
\end{Theorem}
\section{The group of automorphisms of the tree}

For a positive integer $l$, we denote by $X^*$ the set $\{(x_i)_{i=1}^m\mid m\in\mathbb{Z}_{\geq0},x_i\in\mathbb{Z}/l\mathbb{Z}\}$. Connecting $(x_i)_{i=1}^m\in X^*$ with $(x_i)_{i=1}^{m+1}\in X^*$, we can regard $X^*$ as a tree. We denote by $\emptyset$ the root of the tree $X^*$. 

Here after, we only consider the case $l=2$. We denote by $\Aut(X^*)$ the group consisting of automorphisms of $X^*$ satisfying the tree structure of $X^*$. For $w\in X^*$, we set $s^w\in\Aut(X^*)$ as
\begin{align*}
s^w(x)=\begin{cases}
w2&\text{ if }x=w1,\\
w1&\text{ if }x=w2,\\
x&\text{ otherwise}.
\end{cases}
\end{align*}
Naturally, the group $\Aut(X^*)$ is isomorphic to the group generated by $\{s^w\mid w\in X^*\}$ with the relations:
\begin{gather}
(s^w)^2=1\text{ for }w\in X^*,\label{Aut1}\\
s^ws^u=s^{s^w(u)}s^w\text{ for }w,u\in X^*\label{Aut2}.
\end{gather}

For $j=1,2$, let $X^{*}|_j=\{(x_i)_{i=1}^m\mid m\in\mathbb{Z}_{\geq1},x_i\in\mathbb{Z}/2\mathbb{Z},x_1=j\}$. For $\sigma\in\Aut(X^*)$, we set $\sigma|_j\in\Aut(X^*)$ as the restriction of $\sigma$ to $X^{*}|_j$. We also set $\sigma_0\in S_n$ as the restriction of $\sigma$ to $\{1,2\}\subset X^*$. There exists an isomorphism
given by
\[\Aut(X^*)\simeq\mathbb{C}[S_2]\ltimes\Aut(X^*)^{\otimes 2},\qquad\sigma\mapsto(\sigma_0,(\sigma|_1,\sigma|_2)).\]
This relation is called the wreath recursion. 
\begin{Lemma}\label{Lem3.3}
We obtain a homomorphism $\Delta^D\colon\Aut(X^*)\to\Aut(X^*)^{\otimes4}\rtimes\mathbb{C}[S_2]$ by
\[\Delta^D(\sigma)=(\sigma|_{\sigma(1)},\sigma|_{\sigma(2)},\sigma|_{\sigma(2)},\sigma|_{\sigma(1)},\sigma_0),\]
where the action of $\mathbb{C}[S_2]$ on $\Aut(X^*)^{\otimes4}$ can be given by $\sigma_{1,2}(a\otimes b\otimes c\otimes d)=c\otimes d\otimes a\otimes b$.
\end{Lemma}
\begin{proof}
It is enough to show the relation $\Delta^D(\sigma)\Delta^D(\tau)=\Delta^D(\sigma\tau)$. By a direct computation, we obtain
\begin{align*}
\Delta^D(\sigma)\Delta^D(\tau)=(\sigma|_1\tau|_1,\sigma|_2\tau|_2,1)=\Delta^D(\sigma\tau)&\text{ if }\sigma_0=1,\tau_0=1,\\
\Delta^D(\sigma)\Delta^D(\tau)=(\sigma|_1\tau|_2,\sigma|_2\tau|_1,s_{1,2})=\Delta^D(\sigma\tau)&\text{ if }\sigma_0=1,\tau_0=s_{1,2},\\
\Delta^D(\sigma)\Delta^D(\tau)=(\sigma|_2\tau|_2,\sigma|_1\tau|_1,s_{1,2})=\Delta^D(\sigma\tau)&\text{ if }\sigma_0=s_{1,2},\tau_0=1,\\
\Delta^D(\sigma)\Delta^D(\tau)=(\sigma|_2\tau|_1,\sigma|_1\tau|_2,1)=\Delta^D(\sigma\tau)&\text{ if }\sigma_0=s_{1,2},\tau_0=s_{1,2}.
\end{align*}
\end{proof}
In the next section, we will give the Schur-Weyl duality for $\Aut(X^*)$. We prepare one associative algebra.
\begin{Definition}\label{Def907}
For a positive integer $n$, let $G(n)$ be the associative algebra generated by
\[
\{x\otimes\sigma\mid x\in\mathfrak{gl}(n),\sigma\in\Aut(X^*)\}
\]
with the relation:
\begin{equation}
[E_{i,j}\otimes\sigma_1,E_{k,l}\otimes\sigma_2]=\delta_{j,k}E_{i,l}\otimes\sigma_1\sigma_2-\delta_{i,l}E_{k,j}\otimes\sigma_2\sigma_1.\label{def1}
\end{equation}
\end{Definition}
We sometimes denote $E_{i,j}\otimes\sigma$ by $E^\sigma_{i,j}$. The following lemma can be proved by a direct computation.
\begin{Lemma}\label{lem1}
\textup{(1)}\ The homomorphism $U(\mathfrak{gl}(n))$ to $G(n)$ is given by $E_{i,j}\mapsto E^1_{i,j}$.

\textup{(2)}\ There exists a coproduct $\Delta_G$ on $G(n)$ by
\begin{equation}
\Delta_G(E^\sigma_{i,j})=E_{i,j}^\sigma\otimes1+1\otimes E^\sigma_{i,j}\text{ for }1\leq i,j\leq n,\sigma\in \Aut(X^*).
\end{equation} 
\end{Lemma}
\section{Schur-Weyl duality for $\Aut(X^*)$ and $G(m)$}
In this section, for $\sigma=s^{w_1}\cdots s^{w_u}\in\Aut(X^*)$ and $r=1,2$, we set $\sigma_r=s^{rw_1}\cdots s^{rw_u}$. By the definition of $\sigma_r$, we have
\begin{equation}
[\sigma_1,\tau_2]=0,\ 
(\sigma\tau)_r=\sigma_r\tau_r\label{deff1}
\end{equation}
for $\sigma,\tau\in\Aut(X^*)$.
\begin{Theorem}\label{SW-2-1}
For a right module $M$ over $\Aut(X^*)$, there exists an action of $G(n)$ on $M\otimes_{\mathbb{C}[S_2]}(\mathbb{C}^n)^{\otimes 2}$ by
\begin{align*}
E_{i,j}(m\otimes v)&=\sum_{r=1}^2\limits m\otimes E^{(r)}_{i,j}v,\\
E^\sigma_{i,j}(m\otimes v)&=\sum_{r=1}^2\limits m\sigma_r\otimes E^{(r)}_{i,j}v
\end{align*}
for $1\leq i,j\leq n,m\in M,v\in(\mathbb{C}^n)^{\otimes2}$.
\end{Theorem}
\begin{proof}
It is enough to show compatibility with \eqref{def1}.
For $1\leq i,j,k,l\leq n$, $\sigma\in\Aut(X^*)$, $m\in M$ and $v\in(\mathbb{C}^n)^{\otimes 2}$, we obtain
\begin{align*}
&\quad[E_{i,j}^\tau,E_{k,l}^\sigma](m\otimes v)\\
&=m\sigma_1\tau_1\otimes E^{(1)}_{i,j}E^{(1)}_{k,l}v+m\sigma_2\tau_1\otimes E^{(1)}_{i,j}E^{(2)}_{k,l}v+m\sigma_1\tau_2\otimes E^{(2)}_{i,j}E^{(1)}_{k,l}v+m\sigma_2\tau_2\otimes E^{(2)}_{i,j}E^{(2)}_{k,l}v\\
&\quad-m\tau_1\sigma_1\otimes E^{(1)}_{k,l}E^{(1)}_{i,j}v-m\tau_2\sigma_1\otimes E^{(1)}_{k,l}E^{(2)}_{i,j}v-m\tau_1\sigma_2\otimes E^{(2)}_{k,l}E^{(1)}_{i,j}v-m\tau_2\sigma_2\otimes E^{(2)}_{k,l}E^{(2)}_{i,j}v\\
&=\delta_{k,j}m\sigma_1\tau_1\otimes E^{(1)}_{i,l}v+\delta_{k,j}m\sigma_2\tau_2\otimes E^{(2)}_{i,l}v-\delta_{i,l}m\tau_1\sigma_1\otimes E^{(1)}_{k,j}v-\delta_{i,l}m\tau_2\sigma_2\otimes E^{(2)}_{k,j}v\\
&=(\delta_{j,k}E_{i,l}^{\sigma\tau}-\delta_{i,l}E_{k,j}^{\tau\sigma})(m\otimes v),
\end{align*}
where the second and last equalities are due to \eqref{deff1}.
\end{proof}

By Lemma~\ref{lem1} (1), we can regard any module over $G(n)$ as $\mathfrak{gl}(n)$. We say that a $G(n)$-module $N$ is of level $2$ if $N$ is of level $2$ as a $\mathfrak{gl}(n)$-module.

\begin{Theorem}\label{SW128}
Suppose that $n>3$ and $N$ is a $G(n)$-module of level $2$. There exists a right module $M$ over $\Aut(X^*)$ satisfying $M\otimes_{\mathbb{C}[S_2]}(\mathbb{C}^n)^{\otimes2}\simeq N$. 
\end{Theorem}
\begin{proof}
By the definition of $\Aut(X^*)$, there exists a homomorphism from $\mathbb{C}[S_2]$ to $\Aut(X^*)$ by $\sigma_{1,2}\mapsto s^\emptyset$.
Thus, by Theorem~\ref{SW2}, there exists a right module $M$ over $\mathbb{C}[S_2]$ satisfying $M\otimes_{\mathbb{C}[S_2]}(\mathbb{C}^n)^{\otimes2}\simeq N$. The following claim can be proved by the same way as the proof of Lemma~3.14 in \cite{GJU}.
\begin{Claim}\label{Lem2}
For $m\in M,w\in X^*$, there exists $m^{w}\in M$ satisfying the conditions that
\[
E^{s^w}_{i,j}(m\otimes v)=\sum_{r=1}^2m^{rw}\otimes E^{(r)}_{i,j}v\text{ for }1\leq i,j\leq n,v\in(\mathbb{C}^n)^{\otimes2},
\]
and $m^w$ is independent of the choice of $i,j,v$.
\end{Claim}
We set the action of $\Aut(X^*)$ by $ms^w=m^w$ for $w\in X^*\setminus\{\emptyset\}$. We need to prove compatibility with \eqref{Aut1} and \eqref{Aut2}. By the assumption $n>3$, we take $1\leq i<j<k<l\leq n$. 

First, we show compatibility with \eqref{Aut1}. We need to show $m(s^w)^2=m$ for $w\in X^*$ and $m\in M$. In the case $w=\emptyset$, we obtain $m(s^\emptyset)^2=m(\sigma_{1,2})^2=m$.
Suppose that $w\neq\emptyset$.
By a direct computation, we obtain
\begin{align*}
m\otimes (e_{i}\otimes e_{l})=E_{i,k}(m\otimes e_k\otimes e_l)=[E^{w}_{i,j},E^w_{j,k}](m\otimes e_k\otimes e_l)&=ms^{1w}s^{1w}\otimes (e_{i}\otimes e_{l}),
\end{align*}
where the second equality is due to \eqref{def1}.
Then, by Theorem~\ref{SW1} (3), we find the relation $m(s^{1w})^2=m$. Similarly, we can prove $m(s^{2w})^2=m$.

Next, we show the compatibility with \eqref{Aut2}.
We need to show the relation $ms^{w_1}s^{w_2}=ms^{s^{w_1}(w_2)}s^{w_1}$ for $w_1,w_2\in X^*$ and $m\in M$. In the case $w_1=w_2=\emptyset$, it is obvious since the both sides are equal to $m$. In the case $w_1=1w_3,w_2=\emptyset$ for some $w_3\in X^*$, by a direct computation, we obtain
\begin{align*}
&\quad ms^{1w_3}s^\emptyset\otimes e_k\otimes e_i=ms^{1w_3}\otimes e_i\otimes e_k=ms^{1w_3}\otimes E^{(1)}_{i,j}(e_j\otimes e_k)\\
&=E^{w_3}_{i,j}(m\otimes \sigma_{1,2}(e_k\otimes e_j))=E^{w_3}_{i,j}(m\sigma_{1,2}\otimes e_k\otimes e_j)=m\sigma_{1,2}s^{1w_3}\otimes e_k\otimes e_i=ms^\emptyset s^{1w_3}\otimes e_k\otimes e_i.
\end{align*}
Thus, by Theorem~\ref{SW1} (3), we have $ms^{1w_3}s^\emptyset=ms^\emptyset s^{1w_3}$. By the same way, we obtain  
$ms^{2w_3}s^\emptyset=ms^\emptyset s^{2w_3}$.

In the case $w_1=\emptyset,w_2=1w_3$ for some $w_3\in X^*$, by a direct computation, we obtain
\begin{align*}
&\quad ms^\emptyset s^{1w_3}\otimes e_i\otimes e_k=E^{w_3}_{i,j}ms^\emptyset\otimes e_j\otimes e_k=E^{w_3}_{i,j}m\otimes e_k\otimes e_j\\
&=ms^{2w_3}\otimes e_k\otimes e_i=ms^{2w_3}s^\emptyset \otimes e_i\otimes e_k.
\end{align*}
Thus, by Theorem~\ref{SW1} (3), we have $ms^\emptyset s^{1w_3}=m s^{2w_3}s^\emptyset$. By the same way, we obtain  
$ms^\emptyset s^{2w_3}=m s^{1w_3}s^\emptyset$.

In the case $w_1=1w_3,w_2=1w_4$ for some $w_3,w_4\in X^*$, by a direct computation, we obtain
\begin{align*}
&\quad ms^{1w_3}s^{1w_4}\otimes e_{i}\otimes e_{l}=[E^{w_3}_{i,j},E^{w_4}_{j,k}](m\otimes e_{k}\otimes e_{l})\\
&=[E^{s^{w_3}(w_4)}_{i,j},E^{w_3}_{j,k}](m\otimes e_{k}\otimes e_{l})=ms^{1s^{w_3}(w_4)}s^{1w_3}\otimes e_{i}\otimes e_{l},
\end{align*}
where the second equality is due to \eqref{def1}.
Thus, by Theorem~\ref{SW1} (3), we have $ms^{w_1}s^{w_2}=ms^{s^{w_1}(w_2)}s^{w_1}$. The case $w_1=2w_3,w_2=2w_4$ can be proved by the same way.

In the case $w_1=1w_3,w_2=2w_4$ for some $w_3,w_4\in X^*$, by a direct computation, we obtain
\begin{align*}
&\quad ms^{1w_3}s^{2w_4}\otimes e_{i}\otimes e_{j}=E^{w_3}_{i,i}E^{w_4}_{j,j}(m\otimes e_{i}\otimes e_{j})\\
&=E^{w_4}_{j,j}E^{w_3}_{i,i}(m\otimes e_{i}\otimes e_{j})=ms^{2w_4}s^{1w_3}\otimes e_{i}\otimes e_{j},
\end{align*}
where the second equality is due to \eqref{def1}. Thus, by Theorem~\ref{SW1} (3), we have proved the relation $ms^{1w_3}s^{2w_4}=ms^{2w_4}s^{1w_3}$. 
\end{proof}
The bijectivity of the functor 
\[\Hom_{\Aut(X^*)}(M_1,M_2)\to \Hom_{G(m)}(M_1\otimes_{\mathbb{C}[S_2]}(\mathbb{C}^n)^{\otimes2},M_2\otimes_{\mathbb{C}[S_2]}(\mathbb{C}^n)^{\otimes2}),f\mapsto f\otimes\id^{\otimes2}\]
can be proved by the same way as the Schur-Weyl duality for the Yangian and the degenerate affine Hecke algebra. Thus, we obtain the following theorem.
\begin{Theorem}
If $n>3$, the functor 
\[
M\mapsto M\otimes_{\mathbb{C}[S_2]}(\mathbb{C}^n)^{\otimes 2}
\]
gives an equivalence between the category of finite dimensional right modules over $\Aut(X^*)$ to the category of finite dimensional representations of $Y(m)$ can be given by $M\mapsto M\otimes_{\mathbb{C}[S_2]}(\mathbb{C}^n)^{\otimes 2}$. 
\end{Theorem}
\section{The degenerate affine Hecke algebra and the Yangian}
We recall the Yangian and the degenerate affine Hecke algebra shortly. The Yangian is a quantum group, which was a deformation of the current algebra. The Yangian has several presentations. Here, we recall the RTT presentation.
\begin{Definition}
For a complex number $\hbar$, the Yangian $Y_\hbar(\mathfrak{gl}(n))$ is the associative algebra generated by $\{T^{(r)}_{i,j}\mid r\geq0,1\leq i,j\leq n\}$ with the defining relations $T^{(0)}_{i,j}=\delta_{i,j}$ and 
\begin{gather*}
[T^{(r+1)}_{i,j},T^{(s)}_{k,l}]-[T^{(r)}_{i,j},T^{(s+1)}_{k,l}]=\hbar T^{(r)}_{k,j}T^{(s)}_{i,l}-\hbar T^{(s)}_{k,j}T^{(r)}_{i,l}
\end{gather*}
for $r,s\geq0$ and $1\leq i,j,k,l\leq n$.
\end{Definition}

Next, we recall the degenerate affine Hecke algebra of type $A_1$. The degenerate affine Hecke algebra $H_\hbar$ comes from the Cherednik-Dunkl operators, which are operators on $\mathbb{C}[y_1,y_2]$ defined by
\begin{align*}
D_{1}&=y_{1}\dfrac{\partial}{\partial y_{1}}+\hbar+\dfrac{y_{2}}{y_{2}-y_{1}}(1-\theta),\\
D_{2}&=y_{2}\dfrac{\partial}{\partial y_{2}}-\dfrac{y_{2}}{y_{1}-y_{2}}(1-\theta),
\end{align*}
where $\theta$ is an operator on $\mathbb{C}[y_1,y_2]$ which gives the exchange of $y_1$ and $y_2$. 
\begin{Definition}
For a complex number $\hbar$, the associative algebra $H_\hbar$ of type $A_1$ is the associative algebra generated by the symmetric group $S_2$ and $\{u_1,u_2\}$ with the relations:
\begin{gather*}
s_{1,2}u_1-u_2s_{1,2}=-\hbar,\ 
[u_1,u_2]=0.
\end{gather*}
\end{Definition}
We note the element $u_i$ corresponds to $D_{i}$ and $s_{1,2}$ corresponds to $\theta$ (see \cite{Ch}). 

We sometimes use other generators of $H_\hbar$. Let $z_1=u_1$ and $z_2=u_2-\hbar s_{1,2}$. Then, we can rewrite the defining relations as
\begin{equation}
z_1s_{1,2}=s_{1,2}z_2,\ [z_1,z_2]=-\hbar(z_1-z_2)\sigma.
\end{equation}

It was shown that the Yangian $Y(\mathfrak{gl}(n))$ has the Schur-Weyl duality with the degenerate affine Hecke algebra $H_\hbar$ (see \cite{Dri} and \cite{ara}). 
\begin{Theorem}\label{Thm137}
\textup{(1)}\ For a right module $M$ over $H_\hbar$, $SW(M)=M\otimes_{\mathbb{C}[S_l]}(\mathbb{C}^n)^{\otimes l}$ has a $Y(\mathfrak{gl}(n))$-module structure by
\[
T^{(s)}_{i,j}(m\otimes v)=\sum_{r=1}^2\limits m(z_r)^{s-1}\otimes E^{(r)}_{i,j}v
\]
for $1\leq i,j\leq n,r\geq1,m\in M,v\in(\mathbb{C}^n)^{\otimes2}$.

\textup{(2)}\ If $n>3$, the functor $SW$ gives the equivalence between the category of right modules over $H_\hbar$ and the category of finite dimensional representations of the Yangian of level $2$.
\end{Theorem}
\begin{Remark}
Drinfeld \cite{Dri} and Arakawa \cite{ara} gave Theorem~\ref{Thm137} for general $l$.
\end{Remark}
\begin{Lemma}\label{Lem128}
The following relation holds:
\[
z_1^{k+1}z_2^l-z_1^kz_2^{l+1}=\hbar(z_1^kz_2^l-z_1^lz_2^k)\sigma.
\]
\end{Lemma}
\begin{proof}
We can take distinct $1\leq i,j,k,l\leq n$ by the assumption $n\geq4$. By Theorem~\ref{Thm137}, we have an action of $Y(\mathfrak{gl}(n))$ on $H_\hbar\otimes(\mathbb{C}^n)^{\otimes2}$. By a direct computation, we obtain
\begin{align}
([T^{(r+1)}_{i,j},T^{(s)}_{k,l}]-[T^{(r)}_{i,j},T^{(s+1)}_{k,l}])(1\otimes e_j\otimes e_l)
&=([z_1^{r+1},z_2^s]-[z_1^{r+1},z_2^s])\otimes e_i\otimes e_k\label{rel965}
\end{align}
and
\begin{align}
\hbar(T^{(r)}_{k,j}T^{(s)}_{i,l}-T^{(s)}_{k,j}T^{(r)}_{i,l})(1\otimes e_j\otimes e_l)
&=\hbar(z_1^{r}z_2^s-z_1^sz_2^r)\otimes e_k\otimes e_i=\hbar(z_1^{r}z_2^s-z_1^sz_2^r)\sigma\otimes (e_i\otimes e_k)\label{rel966}.
\end{align}
By the definition of the Yangian $Y(\mathfrak{gl}(n))$, the left hand side of $\eqref{rel965}$ and \eqref{rel966}.
Thus, by Theorem~\ref{SW1} (3) and the assumption $i\neq k$, we obtain the relation $z_1^{k+1}z_2^l-z_1^kz_2^{l+1}=\hbar(z_1^kz_2^l-z_1^lz_2^k)\sigma$.
\end{proof}

\section{The degenerate affine Hecke algebra on $X^*$}
In this section, we extend the degenerate affine Hecke algebra $H_\hbar$ to $X^*$. Since the degenerate affine Hecke algebra comes from the Cherednik-Dunkl operators, it is expected that the degenerate affine Hecke algebra associated with $X^*$ relates to operators on $\mathbb{C}[y_{w1},y_{w2}\mid w\in X^*]$ defined by
\begin{align*}
D_{w1}&=y_{w1}\dfrac{\partial}{\partial y_{w1}}+\hbar+\dfrac{y_{w2}}{y_{w2}-y_{w1}}(1-\theta_w),\\
D_{w2}&=y_{w2}\dfrac{\partial}{\partial y_{w2}}-\dfrac{y_{w2}}{y_{w1}-y_{w2}}(1-\theta_w),
\end{align*}
where $\theta_w$ is the exchange of $y_{w1}$ and $y_{w2}$. By a direct computation, we obtain the following relations:
\begin{gather*}
s^wu_x-u_{s^w(x)}s^w=-(\delta_{x,w1}+\delta_{x,w2})\hbar,\qquad 
[u_w,u_x]=0.
\end{gather*}
In particular, we obtain
\[
s^\emptyset u_1-u_2s^\emptyset=-\hbar,\\
s^\emptyset u_w-u_{s^\emptyset(w)}s^\emptyset=0\text{ for }w\neq1,2.
\]
Since this difference is inconvenient for the construction of the Schur-Weyl duality, we change the relations a little. 
\begin{Definition}
For a complex number $\hbar$, we define $\mathbb{H}$ as the associative algebra generated by 
\[\text{Aut}(X^*)\cup\{z_1,z_2\} \cup\{u_w\mid w\in X^*\setminus\{\emptyset,1,2\},i=1,2\}\]
with the relations:
\begin{gather}
s^wz_j-z_{s^w(j)}=0,\label{H1}\\
s^wu_x-u_{w(x)}s^w=(\delta_{x,w1}+\delta_{x,w2})\hbar,\label{H2}\\
[z_1,z_2-\hbar s^\emptyset]=0,\label{H3}\\
[u_w,z_j]=0,\label{H4}\\
[u^{w_1}_{i_1},u^{w_2}_{i_2}]=0.\label{H5}
\end{gather}
\end{Definition}
We will give the Shcur-Weyl duality for $\mathbb{H}$. We prepare some notations.
Let $Q$ be 
\[
\{\mathfrak{p}=(p_1,\cdots,p_m)\mid m\geq0,p_i\in X^*\setminus\{\emptyset\}\}/\sim,
\] 
where $(p_1,\cdots,p_m)\sim(q_1,\cdots,q_u)$ if $m=u$ and there exists $\sigma\in S_m$ satisfying $q_j=p_{\sigma(j)}$ for $1\leq j\leq m$. For $r=1,2$ and $\mathfrak{p}=(p_1,\cdots,p_m)\in Q$, we denote $\mathfrak{p}_r$ by $(rp_1,\cdots,rp_m)$. For $\mathfrak{p}=(p_1,\cdots,p_m)\in Q$, we set $u_{\mathfrak{p}}=u_{p_1}\cdots u_{p_m}\in\mathbb{H}$. 
For $\mathfrak{p_1}\cdots\mathfrak{p_r}\in Q$ and $a_1,\cdots,a_r\in\mathbb{C}$, we set
\[
u_{\sum_{i=1}^r\limits a_i\mathfrak{p_i}}=\sum_{i=1}^r\limits a_iu_{\mathfrak{p_i}}.
\]
For $\mathfrak{p}=(p_1,\cdots,p_m)\in Q$ and $w\in X^*$, we construct $u_{\mathfrak{p}^w}\in\bigoplus_{\mathfrak{q}\in Q}\mathbb{C}\mathfrak{q}$ from $u_{\mathfrak{p}}$ by replacing $u_{w1}$ (resp. $u_{w2}$) with $u_{w2}+\hbar$ (resp. $u_{w1}-\hbar_w$) and $u_{w'}$ with $u_{s^w(w')}$ for $w'\neq w1,w2$ respectively. 
\begin{Definition}
For a positive integer $n$, let $YY(n)$ be the associative algebra generated by
\[
\{T^{\mathfrak{p},\tau,(x)}_{i,j}\mid \mathfrak{p}\in Q\cup\{\emptyset\},\tau\in\Aut(X^*),1\leq i,j\leq n,x\geq1\}\cup\{T^{\mathfrak{p},(x)}_{i,j}\mid \mathfrak{p}\in P_+,1\leq i,j\leq n,x\geq1\}
\]
with the defining relations:
\begin{align}
[T^{\mathfrak{p},\sigma,(1)}_{i,j},T^{\mathfrak{q},\tau,(r)}_{k,l}]&=\delta_{j,k}T^{(\mathfrak{p},\mathfrak{q}^\sigma),\sigma\tau,(r)}_{i,l}-\delta_{i,l}T^{(\mathfrak{p}^\tau,\mathfrak{q}),\tau\sigma,(r)}_{k,j},\label{gat1}\\
[T^{\mathfrak{p},(1)}_{i,j},T^{\mathfrak{q},\tau,(r)}_{k,l}]&=\delta_{j,k}T^{(\mathfrak{p},\mathfrak{q}),\tau,(r)}_{i,l}-\delta_{i,l}T^{(\mathfrak{p}^\tau,\mathfrak{q}),\tau,(r)}_{k,j},\label{gat2}\\
[T^{\mathfrak{p},(1)}_{i,j},T^{\mathfrak{q},(r)}_{k,l}]&=\delta_{j,k}T^{(\mathfrak{p},\mathfrak{q}),\sigma\tau,(r)}_{i,l}-\delta_{i,l}T^{(\mathfrak{p},\mathfrak{q}),\tau\sigma,(r)}_{k,j}\label{gat3}
\end{align}
for $\mathfrak{p}\in P,\sigma\cup\{\emptyset\},\tau\in\Aut(X^*),r\geq1,1\leq i,j,k,l\leq n$ and
\begin{align}
[T^{\mathfrak{p},\sigma,(s+1)}_{i,j},T^{\mathfrak{q},\tau,(r)}_{k,l}]-[T^{\mathfrak{p},\sigma,(s)}_{i,j},T^{\mathfrak{q},\tau,(r+1)}_{k,l}]
&=\hbar T^{\mathfrak{p},\sigma,(s)}_{k,j}T^{\mathfrak{q},\tau,(r)}_{i,l}-\hbar T^{\mathfrak{p},\sigma,(r)}_{k,j}T^{\mathfrak{q},\tau,(s)}_{i,l},\label{gat4}\\
[T^{\mathfrak{p},(s+1)}_{i,j},T^{\mathfrak{q},\tau,(r)}_{k,l}]-[T^{\mathfrak{p},(s)}_{i,j},T^{\mathfrak{q},\tau,(r+1)}_{k,l}]
&=\hbar T^{\mathfrak{p},(s)}_{k,j}T^{\mathfrak{q},\tau,(r)}_{i,l}-\hbar T^{\mathfrak{p},(r)}_{k,j}T^{\mathfrak{q},\tau,(s)}_{i,l},\label{gat5}\\
[T^{\mathfrak{p},(s+1)}_{i,j},T^{\mathfrak{q},(r)}_{k,l}]-[T^{\mathfrak{p},(s)}_{i,j},T^{\mathfrak{q},(r+1)}_{k,l}]
&=\hbar T^{\mathfrak{p},(s)}_{k,j}T^{\mathfrak{q},(r)}_{i,l}-\hbar T^{\mathfrak{p},(r)}_{k,j}T^{\mathfrak{q},(s)}_{i,l}\label{gat6}
\end{align}
for $\mathfrak{p}\in P\cup\{\emptyset\},\tau\in\Aut(X^*),r,s\geq1,1\leq i,j,k,l\leq n$.
\end{Definition}
\begin{Theorem}
For a right module $M$ over $\mathbb{H}$, the action of $YY(n)$ on $M\otimes_{\mathbb{C}[S_2]}(\mathbb{C}^n)^{\otimes2}$ can be defined by
\begin{align*}
T^{\mathfrak{p},(r)}_{i,j}(m\otimes v)&=mu_{\mathfrak{p}_1}z_1^{r-1}\otimes E^{(1)}_{i,j}v+mu_{\mathfrak{p}_2}z_2^{r-1}\otimes E^{(2)}_{i,j}v,\\
T^{\mathfrak{p},\sigma,(r)}_{i,j}(m\otimes v)&=mu_{\mathfrak{p}_1}\sigma_1z_1^{r-1}\otimes E^{(1)}_{i,j}v+mu_{\mathfrak{p}_2}\sigma_2z_2^{r-1}\otimes E^{(2)}_{i,j}v.
\end{align*}
\end{Theorem}
\begin{proof}
It is enough to show the compatibility with \eqref{gat1}-\eqref{gat6}. We only prove the compatibilities with \eqref{gat1} and \eqref{gat4}. 

First, we show the compatibility with \eqref{gat1}. By a direct computation, we obtain
\begin{align*}
[T^{\mathfrak{p},\sigma,(1)}_{i,j},T^{\mathfrak{q},\tau,(r)}_{k,l}](m\otimes v)
&=mu_{\mathfrak{q}_1}\tau_1u_{\mathfrak{p}_1}\sigma_1z_1^{r-1}\otimes E^{(1)}_{i,j}E^{(1)}_{k,l}v+mu_{\mathfrak{q}_2}\tau_2u_{\mathfrak{p}_2}\sigma_2z_2^{r-1}\otimes E^{(2)}_{i,j}E^{(2)}_{k,l}v\\
&\quad-mu_{\mathfrak{p}_1}\sigma_1u_{\mathfrak{q}_1}\tau_1z_1^{r-1}\otimes E^{(1)}_{k,l}E^{(1)}_{i,j}v-mu_{\mathfrak{p}_2}\sigma_2u_{\mathfrak{q}_2}\tau_2z_2^{r-1}\otimes E^{(2)}_{k,l}E^{(2)}_{i,j}v\\
&=mu_{\mathfrak{q}_1}u_{\mathfrak{p}^\tau_1}\tau_1\sigma_1z_1^{r-1}\otimes E^{(1)}_{i,j}E^{(1)}_{k,l}v+mu_{\mathfrak{q}_2}u_{\mathfrak{p}^\tau_2}\tau_2\sigma_2z_2^{r-1}\otimes E^{(2)}_{i,j}E^{(2)}_{k,l}v\\
&\quad-mu_{\mathfrak{p}_1}u_{\mathfrak{q}^\sigma_1}\sigma_1\tau_1z_1^{r-1}\otimes E^{(1)}_{k,l}E^{(1)}_{i,j}v-mu_{\mathfrak{p}_2}u_{\mathfrak{q}^\sigma_2}\sigma_2\tau_2z_2^{r-1}\otimes E^{(2)}_{k,l}E^{(2)}_{i,j}v\\
&=(\delta_{k,j}T^{\mathfrak{p}^\tau,\sigma\tau,(r)}_{i,l}-\delta_{i,l}T^{\mathfrak{q}^\sigma,\tau\sigma,(r)}_{k,j})(m\otimes v),
\end{align*}
where the first equality is due to \eqref{H1}, \eqref{H4} and \eqref{H5} and the second equality is due to \eqref{H2}. Thus, we have proved the compatibility with \eqref{gat1}.

Next, we show the compatibility with \eqref{gat4}. By a direct computation, we obtain
\begin{align*}
&\quad([T^{\mathfrak{p},\sigma,(r+1)}_{i,j},T^{\mathfrak{q},\tau,(s)}_{k,l}]-[T^{\mathfrak{p},\sigma,(r)}_{i,j},T^{\mathfrak{q},\tau,(s+1)}_{k,l}])(m\otimes v)\\
&=mu_{\mathfrak{q}_1}\tau_1u_{\mathfrak{p}_2}\sigma_2(z_1^sz_2^{r+1}-z_1^{s+1}z_2^r)\otimes E^{(2)}_{i,j}E^{(1)}_{k,l}v\\
&\quad+mu_{\mathfrak{q}_2}\tau_2u_{\mathfrak{p}_1}\sigma_1(z_2^sz_1^{r+1}-z_2^{s+1}z_1^r)\otimes E^{(1)}_{i,j}E^{(2)}_{k,l}v\\
&=mu_{\mathfrak{q}_1}\tau_1u_{\mathfrak{p}_2}\sigma_2(z_1^sz_2^r-z_1^rz_2^s)\sigma\otimes E^{(2)}_{k,j}E^{(1)}_{i,l}v\\
&\quad+mu_{\mathfrak{q}_2}\tau_2u_{\mathfrak{p}_1}\sigma_1((z_2)^{s}z_1^r-z_2^rz_1^{s})\otimes E^{(1)}_{i,j}E^{(2)}_{k,l}v\\
&=(T^{\mathfrak{p},\sigma,(r)}_{k,j}T^{\mathfrak{q},\tau,(s)}_{i,l}-T^{\mathfrak{p},\sigma,(s)}_{k,j}T^{\mathfrak{q},\tau,(r)}_{i,l})(m\otimes v),
\end{align*}
where the first equality is due to \eqref{H1}, \eqref{H4}, \eqref{H5} and Lemma~\ref{Lem128} and the second equality is due to \eqref{H2}. Thus, we have proved the compatibility with \eqref{gat4}.
\end{proof}
There exists a homomorphism from $U(\mathfrak{gl}(n))$ to $YY(n)$ given by $E_{i,j}\mapsto T^{\emptyset,(1)}_{i,j}$. We say that a $YY(n)$-module $M$ is of type $2$ if $M$ is of type $2$ as a $YY(n)$-module.
\begin{Theorem}
If $n>3$, for a right module $N$ over $YY(n)$ of type $2$, there exists an $\mathbb{H}$-right module $M$ satisfying $M\otimes_{\mathbb{C}[S_2]}(\mathbb{C}^n)^{\otimes2}\simeq N$.
\end{Theorem}
\begin{proof}
There exists a homomorphism from $\Aut(X^*)$ to $\mathbb{H}$ by $\sigma\mapsto \sigma$.
By Theorem~\ref{SW128}, there exists a right $\mathbb{C}[S_2]$-module $M$ over $\Aut(X^*)$ satisfying $M\otimes_{\mathbb{C}[S_2]}(\mathbb{C}^n)^{\otimes2}\simeq N$. By the same way as the proof of Lemma~3.14 in \cite{GJU}, for $m\in M$ and $w\in X^*\setminus\{\emptyset\}$, we find that there exist $m^{w}\in M$ satisfying conditions that
\begin{align*}
T^{(w),(2)}_{i,j}(m\otimes v)&=\sum_{a=1}^2\limits m_a\otimes E^{(a)}_{i,j}v\text{ for }1\leq i,j\leq n,v\in(\mathbb{C}^n)^{\otimes2},\\
T^{(w),(1)}_{i,j}(m\otimes v)&=\sum_{a=1}^2\limits m^{w}_a\otimes E^{(a)}_{i,j}v\text{ for }1\leq i,j\leq n,v\in(\mathbb{C}^n)^{\otimes2}
\end{align*}
and $m^w$ is independent of the choice of $v$. We set $mu_{rw}=m^{w}_r$ and $mz_r=m_r$ for $w\in X^*\setminus\{\emptyset\}$. It is enough to show that this action is compatible with \eqref{H1}-\eqref{H5}.
We only show the compatibility with the relations $[u_{1w_1},u_{2w_2}]=0$ and $s_{1w_1}u_{1w_2}=u_{1s^{w_1}(w_2)}s_{1w_1}+\delta_{w_2,w_11}\hbar s^{1w_1}-\delta_{w_2,w_12}\hbar s^{1w_1}$. The other compatibilities can be proved by the similar way. By the assumption $n>3$, we can take $1\leq i\leq j\leq k\leq l\leq n$.

By a direct computation, we obtain
\begin{align*}
&\quad mu_{1w_1}u_{2w_2}\otimes e_i\otimes e_k=mu_{1w_1}u_{2w_2}\otimes E^{(1)}_{i,j}E^{(2)}_{k,l}e_j\otimes e_l\\
&=T^{(w_2),(1)}_{k,l}T^{(w_1),(1)}_{i,j}m\otimes e_j\otimes e_l\\
&=T^{(w_1),(1)}_{i,j}T^{(w_2),(1)}_{k,l}m\otimes e_j\otimes e_l\\
&=mu_{2w_2}u_{1w_1}\otimes e_i\otimes e_k,
\end{align*}
where the third equality is due to
\begin{align*}
[T^{(w_2),(1)}_{k,l},T^{(w_1),(1)}_{i,j}]&=\delta_{i,l}T^{(w_1,w_2),(1)}_{k,j}-\delta_{j,k}T^{(w_1,w_2),(1)}_{i,l}=0.
\end{align*}
By Theorem~\ref{SW1} (3), we have proved the compatibility with $[u_{1w_1},u_{2w_2}]=0$.

By a direct computation, we obtain
\begin{align*}
&\quad ms_{1w_1}u_{1w_2}\otimes e_i\otimes e_l=[T^{(w_2),(1)}_{i,j},T^{\emptyset,s^{w_1},(1)}_{j,k}]m\otimes e_k\otimes e_l\\
&=(T^{(w_1,s^{w_1}(w_2)),s^{w_1},(1)}_{i,k}+\hbar\delta_{w_11,w_2}T^{s^{w_1},(1)}_{i,k}-\hbar\delta_{w_12,w_2}T^{s^{w_1},(1)}_{i,k})(m\otimes e_k\otimes e_l)\\
&=m(u_{1s^{w_1}(w_2)}s_{1w_1}+\delta_{w_2,w_11}\hbar s^{1w_1}-\delta_{w_2,w_12}\hbar s^{1w_1})\otimes e_i\otimes e_l.
\end{align*}
By Theorem~\ref{SW1} (3), we have proved the compatibility with $[u_{1w_1},u_{2w_2}]=0$ and $s_{1w_1}u_{1w_2}=u_{1s^{w_1}(w_2)}s_{1w_1}+\delta_{w_2,w_11}\hbar s^{1w_1}-\delta_{w_2,w_12}\hbar s^{1w_1}$.
\end{proof}
Similarly to the wreath recursion of $\Aut(X^*)$, we expect that $YY(n)$ has a homomorphism corresponding to the wreath recursion.
The Yangian $Y_\hbar(\mathfrak{gl}(n))$ has a Hopf algebra structure whose coproduct 
\[
\Delta\colon Y(\mathfrak{gl}(n))\to Y(\mathfrak{gl}(n))\otimes Y(\mathfrak{gl}(n))
\]
can be given by
\[
\Delta(T^{(r)}_{i,j})=T^{(r)}_{i,j}\otimes 1+1\otimes T^{(r)}_{i,j}+\hbar\sum_{u=0}^r\limits\sum_{x=1}^n\limits T^{(u-r)}_{i,x}\otimes T^{(r)}_{x,j}.
\]
Thus, we expect that we can construct a homomorphism for $YY(n)$ corresponding to the wreath recursion by the similar way to the coproduct for the original Yangian. However, it is difficult. From the next section, we change the definition of $YY(n)$ a little and show that it has the homomorphism corresponding the coproduct.

\section{A different product on $\Aut(X^*)$}
In this section, we give another group structure of $\Aut(X^*)$.
Let us set the product on $\Aut(X^*)$ by
\[
\sigma\ast\tau=\tau_0\sigma\tau_0\tau.
\]
Since we obtain
\begin{gather*}
(\sigma\ast\tau)\ast\omega=\omega_0(\tau_0\sigma\tau_0\tau)\omega_0\omega=\omega_0\tau_0\sigma\tau_0\omega_0(\omega_0\tau\omega_0\omega)=\sigma\ast(\tau\ast\omega),\\
1\ast\sigma=\sigma_01\sigma_0\sigma=\sigma=\sigma\ast1,\\
\sigma\ast(\sigma_0\sigma^{-1}\sigma_0)=\sigma_0\sigma\sigma_0(\sigma_0\sigma^{-1}\sigma_0)=1,\\
(\sigma_0\sigma^{-1}\sigma_0)\ast\sigma=\sigma_0(\sigma_0\sigma^{-1}\sigma_0)\sigma_0\sigma=1,
\end{gather*}
by a direct computation, we find that $\Aut(X^*)$ becomes a group with the product $\ast$. 
We denote this group by $\Aut(X^*)^\ast$. 

By the definition of the product $\ast$, we obtain the following lemma.
\begin{Lemma}\label{Lem980}
\textup{(1)}\ The relation $(\sigma\ast\tau)|_r=\sigma|_r\tau|_r$ holds.

\textup{(2)}\ There exists a homomorphism
\[
\widehat{\Delta}\colon \Aut(X^*)^\ast\to\Aut(X^*)^{\otimes2}\otimes\mathbb{C}[S_2],\ \sigma\mapsto\sigma|_1\otimes\sigma|_2\otimes\sigma_0.
\]

\textup{(3)}\ There exists a homomorphism
\[
\Phi\colon\Aut(X^*)^\ast\to\mathbb{C}[S_2]\ltimes\Aut(X^*),\ \sigma\mapsto(\sigma_0,\sigma),
\]
where the action of $\mathbb{C}[S_2]$ on $\Aut(X^*)^\ast$ is given by
\[
\sigma_{1,2}\colon \tau=\sigma_{1,2}\tau\sigma_{1,2}\text{ for }tau\in\Aut(X^*)^\ast.
\]
\end{Lemma}
\begin{proof}
\textup{(1)}\ By the definition, we have 
\[
(\sigma\ast\tau)|_r=\tau_0|_{\sigma(r)}\sigma|_r\tau_0|_{\tau(r)}\tau|_r=\sigma|_r\tau|_r
\]
for $\sigma,\tau\in\Aut(X^*)$.

\textup{(2)}\ By (1), we obtain
\[
\widehat{\Delta}(\sigma)\widehat{\Delta}(\tau)=\sigma|_1\tau|_1\otimes\sigma|_2\tau|_2\otimes\sigma_0\tau_0=(\sigma\ast\tau)|_1\otimes(\sigma\ast\tau)|_2\otimes\sigma_0\tau_0=\widehat{\Delta}(\sigma\ast\tau)
\]
for $\sigma,\tau\in\Aut(X^*)$.

\textup{(3)}\ By the definition, we obtain
\begin{align*}
&\quad \Phi(\sigma)\Phi(\tau)=((\tau_0\sigma\tau_0)\ast\tau,\sigma_0\tau_0)=(\sigma\tau,\sigma_0\tau_0)=\Psi(\sigma\tau)
\end{align*}
for $\sigma,\tau\in\Aut(X^*)$.
\end{proof}
Similarly to Definition~\ref{Def907}, we can define the associative algebra generated by $\{E^\sigma_{i,j}\mid\sigma\in\Aut(X^*),1\leq i,j\leq n\}$ with the relation:
\begin{equation}
[E_{i,j}^\sigma,E_{k,l}^\tau]=\delta_{j,k}E_{i,l}^{\sigma\ast\tau}-\delta_{i,l}E_{k,j}^{\tau\ast\sigma}.\label{def2}    
\end{equation}
We denote this algebra by $GG(n)$. 
In the rest of this section, we will give the Schur-Weyl duality for $GG(n)$.
\begin{Definition}
Let $A$ be the associative algebra generated by $S_2\cup\{\theta^r_w,\omega^r_w,\mid r=1,2,w\in X^*\}\}\cup\{\delta^1,\delta^2\}$ with the relations:
\begin{gather}
 [\theta^1_w,\theta^2_x]=0,\ [\theta^r_w,\omega^s_x]=0,\ [\omega^1_w,\omega^2_x]=0,\label{A1}\\
 (\theta^r_w)^2=1,\ (\omega^r_w)^2=1,\label{A2}\\
 \theta^r_{w_1}\theta^r_{w_2}=\theta^r_{s^{w_2}(w_1)}\theta^r_{w_1},\label{A3}\\
 \omega^r_{w_1}\omega^r_{w_2}=\omega^r_{s^{w_2}(w_1)}\omega^r_{w_1},\label{A3-1}\\
 [\delta^r,\theta^v_w]=0,\ [\delta^r,\theta^v_w]=0,[\delta^1,\delta^2]=0,\label{A4}\\
 s_{1,2}\theta^r_{w_1}=\theta^{s_{1,2}(r)}_{w_1}s_{1,2},\  s_{1,2}\omega^r_{w_1}=\omega^{s_{1,2}(r)}_{w_1}s_{1,2},\  s_{1,2}\delta^r=\delta^{s_{1,2}(r)}s_{1,2}.\label{A5}
\end{gather}
\end{Definition}
\begin{Remark}
From the definition, $A$ is isomorphic to $\mathbb{C}[S_2]\ltimes(\mathbb{C}[z_1,z_2]\otimes\mathbb{C}[\Aut(X^*)])$.
\end{Remark}
For $\sigma=\sigma_0\sigma|_1\sigma|_2\in\Aut(X^*)$ and $r=1,2$, we set $\theta^r_{\sigma}=\theta^r_{w_1}\cdots\theta^r_{w_m}(\theta^r_\emptyset)^{\sigma_0}$ and $\omega^r_{\sigma}=\omega^r_{w_1}\cdots\omega^r_{w_m}(\theta^r_\emptyset)^{\sigma_0}$ if $\sigma|_1=s^{w_1}\cdots s^{w_m}$ and $\sigma|_2=s^{x_1}\cdots s^{x_m}$. By the definition, we have
\begin{gather}
\theta^r_\sigma\theta^r_\tau=\theta^r_{\tau_0\sigma\tau_0\tau},\label{grorious1}\\
\omega^r_\sigma\omega^r_\tau=\omega^r_{\tau\sigma\tau_0\tau}.\label{grorious2}
\end{gather}
\begin{Theorem}\label{SW-2-2}
For a right module $M$ over $A$, there exists an action of $GG(n)$ on $M\otimes_{\mathbb{C}[S_2]}(\mathbb{C}^n)^{\otimes 2}$ by
\begin{align*}
E^1_{i,j}(m\otimes v)&=\sum_{r=1}^2\limits m\otimes E^{(r)}_{i,j}v,\\
E^\sigma_{i,j}(m\otimes v)&=\sum_{r=1}^2\limits m\theta^r_{\sigma}\omega^r_{\sigma}(\delta^r)^{\sigma_0}\otimes E^{(r)}_{i,j}v\text{ for }\sigma\neq\emptyset
\end{align*}
for $1\leq i,j\leq n,m\in M,v\in(\mathbb{C}^n)^{\otimes2}$.
\end{Theorem}
\begin{proof}
It is enough to show compatibility with \eqref{def2}. We only show the case $\sigma,\tau\neq\emptyset$. By a direct computation, we obtain
\begin{align*}
[E_{i,j}^\tau,E_{k,l}^\sigma](m\otimes v)
&=m\theta^1_\sigma\omega^1_\sigma\theta^1_\tau\omega^1_\tau(\delta^1)^{\sigma_0+\tau_0}\otimes E^{(1)}_{i,j}E^{(1)}_{k,l}v+m\theta^1_\sigma\omega^1_\sigma\theta^2_\tau\omega^2_\tau(\delta^1)^{\sigma_0}(\delta^2)^{\tau_0}\otimes E^{(2)}_{i,j}E^{(1)}_{k,l}v\\
&\quad+m\theta^2_\sigma\omega^2_\sigma\theta^1_\tau\omega^1_\tau(\delta^1)^{\sigma_0}(\delta^2)^{\tau_0}\otimes E^{(1)}_{i,j}E^{(2)}_{k,l}v+m\theta^2_\sigma\omega^2_\sigma\theta^2_\tau\omega^2_\tau(\delta^2)^{\sigma_0+\tau_0}\otimes E^{(2)}_{i,j}E^{(2)}_{k,l}v\\
&\quad-m\theta^1_\sigma\omega^1_\sigma\theta^1_\tau\omega^1_\tau(\delta^1)^{\sigma_0+\tau_0}\otimes E^{(1)}_{k,l}E^{(1)}_{i,j}v-m\theta^1_\sigma\omega^1_\sigma\theta^2_\tau\omega^2_\tau(\delta^2)^{\sigma_0}(\delta^1)^{\tau_0}\otimes E^{(1)}_{k,l}vE^{(2)}_{i,j}\\
&\quad-m\theta^2_\sigma\omega^2_\sigma\theta^1_\tau\omega^1_\tau(\delta^1)^{\sigma_0}(\delta^2)^{\tau_0}\otimes E^{(2)}_{k,l}vE^{(1)}_{i,j}-m\theta^2_\sigma\omega^2_\sigma\theta^2_\tau\omega^2_\tau(\delta^2)^{\sigma_0+\tau_0}\otimes E^{(2)}_{k,l}E^{(2)}_{i,j}v\\
&=\delta_{k,j}m\theta^1_\sigma\omega^1_\sigma\theta^1_\tau\omega^1_\tau(\delta^1)^{\sigma_0+\tau_0}\otimes E^{(1)}_{i,l}v+\delta_{k,j}m\theta^2_\sigma\omega^2_\sigma\theta^2_\tau\omega^2_\tau(\delta^2_\emptyset)^{\sigma_0+\tau_0}\otimes E^{(2)}_{i,l}v\\
&\quad-\delta_{i,l}m\theta^1_\tau\omega^1_\tau\theta^1_\sigma\omega^1_\sigma(\theta^1_\emptyset)^{\sigma_0+\tau_0}\otimes E^{(1)}_{k,j}v-\delta_{i,l}m\theta^2_\tau\omega^2_\tau\theta^2_\sigma\omega^2_\sigma(\theta^2_\emptyset)^{\sigma_0+\tau_0}\otimes E^{(2)}_{k,j}v\\
&=(\delta_{j,k}E_{i,l}^{\tau_0\sigma\tau_0\tau}-\delta_{i,l}E_{k,j}^{\sigma_0\tau\sigma_0\sigma})(m\otimes v),
\end{align*}
where the second equality is due to \eqref{A1},\eqref{A3} and \eqref{A4} and the last equality is due to \eqref{grorious1} and \eqref{grorious2}.
\end{proof}
There exists a homomorphism from $U(\mathfrak{gl}(n))$ to $GG(n)$ given by $E_{i,j}\mapsto E^1_{i,j}$. We say that $M$ is a $GG(n)$-module $N$ of type $2$ if $M$ is of type $2$ as a $\mathfrak{gl}(n)$-module.
\begin{Theorem}
If $n>3$, for a $GG(n)$-finite dimensional module $N$ of type $2$, there exists an $A$-right module $M$ satisfying $M\otimes_{\mathbb{C}[S_2]}(\mathbb{C}^n)^{\otimes2}$.
\end{Theorem}
\begin{proof}
By Theorem~\ref{SW1} (3), there exists a $\mathbb{C}[S_2]$-right module $M$ satisfying $M\otimes_{\mathbb{C}[S_2]}(\mathbb{C}^n)^{\otimes2}$. By the same way as the proof of Lemma~3.14 in \cite{GJU}, we find the following claim. 
\begin{Claim}
For $m\in M,w\in X^*$, there exists $m^{w}\in M$ satisfying the conditions that
\begin{align*}
E^\emptyset_{i,j}(m\otimes v)&=\sum_{r=1}^2m^r_{\emptyset}\otimes E^{(r)}_{i,j}v\text{ for }1\leq i,j\leq n,v\in(\mathbb{C}^n)^{\otimes2},\\
E^{s^{1w}}_{i,j}(m\otimes v)&=\sum_{r=1}^2m^r_{1w}\otimes E^{(r)}_{i,j}v\text{ for }1\leq i,j\leq n,v\in(\mathbb{C}^n)^{\otimes2},\\
E^{s^{2w}}_{i,j}(m\otimes v)&=\sum_{r=1}^2m^r_{2w}\otimes E^{(r)}_{i,j}v\text{ for }1\leq i,j\leq n,v\in(\mathbb{C}^n)^{\otimes2}
\end{align*}
and $m^r_w$ is independent of the choice of $i,j,v$. We set $m\delta^r=m^r_\emptyset$, $m\theta^r_w=m^r_{1w}$ and $m\omega^r_w=m^r_{2w}$. 
\end{Claim}
We set $m\delta^r=m^r_\emptyset$, $m\theta^r_w=m^r_{1w}$ and $m\omega^r_w=m^r_{2w}$. Then, it is enough to show the compatibility with \eqref{A1}-\eqref{A5}.

Here, we only show the compatibility with $[\theta^1_w,\omega^1_x]=0$. The other relations can be proved by the same way. By the assumption $n>3$, we can take $1\leq i<j<k<l\leq n$. By a direct computation, for $m\in M$, we obtain
\begin{align*}
&m\theta^1_w\omega^1_x\otimes e_i\otimes e_l=m\theta^1_w\omega^1_x\otimes E^{(1)}_{i,j}E^{(1)}_{j,k}(e_k\otimes e_l)=[E^{s^{1w}}_{i,j},E^{s^{2x}}_{j,k}](m\otimes e_k\otimes e_l)\\
&\qquad=[E^{s^{2x}}_{i,j},E^{s^{1w}}_{j,k}](m\otimes e_k\otimes e_l)=m\omega^1_x\theta^1_w\otimes E^{(1)}_{i,j}E^{(1)}_{j,k}(e_k\otimes e_l)=m\omega^1_x\theta^1_w\otimes e_i\otimes e_l,
\end{align*}
where the third equality is due to $[E^{s^{1w}}_{i,j},E^{s^{2x}}_{j,k}]=E^{s^{1w}s^{2x}}_{i,k}=E^{s^{2x}s^{1w}}_{i,k}=[E^{s^{2x}}_{i,j},E^{s^{1w}}_{j,k}]$. Thus, we have proved the compatibility with $[\theta^1_w,\omega^1_x]=0$.
\end{proof}

\section{Extension of the Yangian}
Now, we give the definition of a candidate of extensions of the Yangian. Let $X=X^*\sqcup\{\cdot\}$ and $P=\{(a_w)_{w\in X}\mid \sum_{w\in X}a_w<\infty\}$. 
We denote by $(a_w)_x$ the $x$-th term of $(a_w)$. We can define an action of $\Aut(X^*)$ on $P$ by $(\sigma(a_w))_\cdot=(a_w)_\cdot$ and $(\sigma(a_w))_x=(a_w)_{\sigma(x)}$ for $x\in X$ and $\sigma\in\Aut(X^*)$. We set $(a_w)|_1$ and $(a_w)|_2$ as
\begin{align*}
((a_w)|_1)_x&=\begin{cases}
a_{1m}&\text{ if }x\neq\cdot,\\
\sum_{y\not\in X^*|_1}a_y&\text{ if }x=\cdot,
\end{cases}\qquad((a_w)|_1)_x=\begin{cases}
a_{2m}&\text{ if }x\neq\cdot,\\
\sum_{y\not\in X^*|_2}a_y&\text{ if }x=\cdot.
\end{cases}
\end{align*}
Let $(0)\in P$ and $(\delta_w)\in P$ be $(0)_w=0$ and $(\delta_w)_x=\delta_{x,w}$ respectively.
\begin{Definition}
For a complex number $\hbar$, let $Y_\hbar(n)$ be the associative algebra generated by
\[
\{T^{\mathfrak{p}}_{i,j}\mid1\leq i,j\leq n,\mathfrak{p}\in P\}
\]
with the relations:
\begin{align}
[T^{\mathfrak{p}+(\delta_w)}_{i,j},T^{\mathfrak{q}}_{k,l}]-[T^{\mathfrak{p}}_{i,j},T^{\mathfrak{q}+(\delta_w)}_{k,l}]
&=\hbar(T^{\mathfrak{p}}_{k,j}T^{\mathfrak{q}}_{i,l}-T^{\mathfrak{q}}_{k,j}T^{\mathfrak{p}}_{i,l}),\label{gatt1}\\
[T^{\delta_{\cdot}}_{i,j},T^{\mathfrak{p}}_{k,l}]
&=\delta_{j,k}T^{\mathfrak{p}+\delta{\cdot}}_{i,l}-\delta_{i,l}T^{\mathfrak{p}+\delta{\cdot}}_{i,l},\label{gatt2}\\
T^{(0)}_{i,j}&=\delta_{i,j}.\label{gatt3}
\end{align}
\end{Definition}
\begin{Theorem}\label{Thm1}
There exists a homomorphism 
\[\Delta\colon Y_\hbar(n)\to Y_\hbar(n)\otimes Y_\hbar(n)\]
given by
\begin{align*}
\Delta(T^{\mathfrak{p}}_{i,j})&=T^{\mathfrak{p}|_1}_{i,j}\otimes1+1\otimes T^{\mathfrak{p}|_2}_{i,j}+\hbar\sum_{\substack{\mathfrak{r}+\mathfrak{s}=\mathfrak{p},\\\mathfrak{r},\mathfrak{s}\neq(0)}}T^{\mathfrak{r}|_1}_{i,x}\otimes T^{\mathfrak{s}|_2}_{x,j}.
\end{align*}
\end{Theorem}
\begin{proof}
It is enough to show the compatibility with \eqref{gatt1}-\eqref{gatt3}. The compatibility with \eqref{gatt1} is trivial. As for the compatibility with \eqref{gatt1}, by a direct computation, we obtain
\begin{align*}
&\quad[\Delta(T^{\mathfrak{p}+(\delta_w)}_{i,j}),\Delta(T^{\mathfrak{q}}_{k,l})]-[\Delta(T^{\mathfrak{p}}_{i,j}),\Delta(T^{\mathfrak{q}+(\delta_w)}_{k,l})]\\
&=([T^{\mathfrak{p}+(\delta_w)|_1}_{i,j},T^{\mathfrak{q}|_1}_{k,l}]-[T^{\mathfrak{p}|_1}_{i,j},T^{\mathfrak{q}+(\delta_w)|_1}_{k,l}])\otimes1+1\otimes ([T^{\mathfrak{p}+(\delta_w)|_2}_{i,j},T^{\mathfrak{q}|_2}_{k,l}]-[T^{\mathfrak{p}|_2}_{i,j},T^{\mathfrak{q}+(\delta_w)|_2}_{k,l}])\\
&\quad+\hbar\sum_{y=1}^n\limits\sum_{\substack{\mathfrak{r}+\mathfrak{s}=\mathfrak{q},\\\mathfrak{r},\mathfrak{s}\neq(0)}}([T^{\mathfrak{p}+(\delta_w)|_1}_{i,j},T^{\mathfrak{r}|_1}_{k,y}]-[T^{\mathfrak{p}|_1}_{i,j},T^{\mathfrak{r}+(\delta_w)|_1}_{k,y}])\otimes T^{\mathfrak{s}|_2}_{y,l}\\
&\quad+\hbar\sum_{y=1}^n\limits\sum_{\substack{\mathfrak{r}+\mathfrak{s}=\mathfrak{q},\\\mathfrak{r},\mathfrak{s}\neq(0)}}T^{\mathfrak{r}|_1}_{k,y}\otimes ([T^{\mathfrak{p}+(\delta_w)|_2}_{i,j},T^{\mathfrak{s}|_2}_{y,l}]-[T^{\mathfrak{p}|_2}_{i,j},T^{\mathfrak{s}+(\delta_w)|_2}_{y,l}])\\
&\quad+\hbar\sum_{x=1}^n\limits\sum_{\substack{\mathfrak{r}+\mathfrak{s}=\mathfrak{p},\\\mathfrak{r},\mathfrak{s}\neq(0)}}([T^{\mathfrak{r}+(\delta_w)|_1}_{i,x},T^{\mathfrak{q}|_1}_{k,l}]-[T^{\mathfrak{r}|_1}_{i,x},T^{\mathfrak{q}+(\delta_w)|_1}_{k,l}])\otimes T^{\mathfrak{s}|_2}_{x,j}\\
&\quad+\hbar\sum_{x=1}^n\limits\sum_{\substack{\mathfrak{r}+\mathfrak{s}=\mathfrak{p},\\\mathfrak{r},\mathfrak{s}\neq(0)}} T^{\mathfrak{r}|_1}_{i,x}\otimes ([T^{\mathfrak{p}+(\delta_w)|_2}_{x,j},T^{\mathfrak{q}|_2}_{k,l}]-[T^{\mathfrak{p}|_2}_{x,j},T^{\mathfrak{q}+(\delta_w)|_2}_{k,l}])\\
&\quad+\hbar^2\sum_{x,y=1}^n\limits\sum_{\substack{\mathfrak{r}+\mathfrak{s}=\mathfrak{p},\\\mathfrak{r},\mathfrak{s}\neq(0)}}\sum_{\substack{\mathfrak{a}+\mathfrak{b}=\mathfrak{q}}}\limits([T^{\mathfrak{r}+\delta_w|_1}_{i,x},T^{\mathfrak{s}|_1}_{k,y}]-[T^{\mathfrak{r}|_1}_{i,x},T^{\mathfrak{s}+\delta_w|_1}_{k,y}])\otimes T^{\mathfrak{s}|_2}_{x,j}T^{\mathfrak{b}|_2}_{y,l}\\
&\quad+\hbar^2\sum_{x,y=1}^n\limits\sum_{\substack{\mathfrak{r}+\mathfrak{s}=\mathfrak{p},\\\mathfrak{r},\mathfrak{s}\neq(0)}}\sum_{\substack{\mathfrak{a}+\mathfrak{b}=\mathfrak{q},\\\mathfrak{a},\mathfrak{b}\neq(0)}}\limits T^{\mathfrak{a}|_1}_{k,y}T^{\mathfrak{r}|_1}_{i,x}\otimes ([T^{\mathfrak{s}+(\delta_w)|_2}_{x,j},T^{\mathfrak{b}|_2}_{y,l}]-[T^{\mathfrak{s}|_2}_{x,j},T^{\mathfrak{b}+(\delta_w)|_2}_{y,l}])\\
&=(T^{\mathfrak{p}|_1}_{k,j}T^{\mathfrak{q}|_1}_{i,l}-T^{\mathfrak{q}|_1}_{k,j}T^{\mathfrak{p}|_1}_{i,l})\otimes1+1\otimes (T^{\mathfrak{p}|_2}_{k,j}T^{\mathfrak{q}|_2}_{i,l}-T^{\mathfrak{q}|_2}_{k,j}T^{\mathfrak{p}|_2}_{i,l})\\
&\quad+\hbar\sum_{y=1}^n\limits\sum_{\substack{\mathfrak{r}+\mathfrak{s}=\mathfrak{q},\\\mathfrak{r},\mathfrak{s}\neq(0)}}(T^{\mathfrak{p}|_1}_{k,j}T^{\mathfrak{r}|_1}_{i,y}-T^{\mathfrak{r}|_1}_{k,j}T^{\mathfrak{p}|_1}_{i,y})\otimes T^{\mathfrak{s}|_2}_{y,l}+\hbar\sum_{y=1}^n\limits\sum_{\substack{\mathfrak{r}+\mathfrak{s}=\mathfrak{q},\\\mathfrak{r},\mathfrak{s}\neq(0)}}T^{\mathfrak{r}|_1}_{k,y}\otimes (T^{\mathfrak{p}|_2}_{k,j}T^{\mathfrak{s}|_2}_{y,l}-T^{\mathfrak{s}|_2}_{k,j}T^{\mathfrak{p}|_2}_{y,l})\\
&\quad+\hbar\sum_{x=1}^n\limits\sum_{\substack{\mathfrak{r}+\mathfrak{s}=\mathfrak{p},\\\mathfrak{r},\mathfrak{s}\neq(0)}}(T^{\mathfrak{r}|_1}_{k,x}T^{\mathfrak{q}|_1}_{i,l}-T^{\mathfrak{q}|_1}_{k,x}T^{\mathfrak{r}|_1}_{i,l})\otimes T^{\mathfrak{s}|_2}_{x,j}+\hbar\sum_{x=1}^n\limits\sum_{\substack{\mathfrak{r}+\mathfrak{s}=\mathfrak{p},\\\mathfrak{r},\mathfrak{s}\neq(0)}}T^{\mathfrak{r}|_1}_{i,x}\otimes (T^{\mathfrak{p}|_2}_{k,j}T^{\mathfrak{q}|_2}_{x,l}-T^{\mathfrak{q}|_2}_{k,j}T^{\mathfrak{p}|_2}_{x,l})\\
&\quad+\hbar^2\sum_{x,y=1}^n\limits\sum_{\substack{\mathfrak{r}+\mathfrak{s}=\mathfrak{p},\\\mathfrak{r},\mathfrak{s}\neq(0)}}\sum_{\substack{\mathfrak{a}+\mathfrak{b}=\mathfrak{q},\\\mathfrak{a},\mathfrak{b}\neq(0)}}\limits(T^{\mathfrak{r}|_1}_{k,x}T^{\mathfrak{s}|_1}_{i,y}-T^{\mathfrak{s}|_1}_{k,x}T^{\mathfrak{r}|_1}_{i,y})\otimes T^{\mathfrak{s}|_2}_{x,j}T^{\mathfrak{b}|_2}_{y,l}\\
&\quad+\hbar^2\sum_{x,y=1}^n\limits\sum_{\substack{\mathfrak{r}+\mathfrak{s}=\mathfrak{p},\\\mathfrak{r},\mathfrak{s}\neq(0)}}\sum_{\substack{\mathfrak{a}+\mathfrak{b}=\mathfrak{q},\\\mathfrak{a},\mathfrak{b}\neq(0)}}\limits T^{\mathfrak{a}|_1}_{k,y}T^{\mathfrak{r}|_1}_{i,x}\otimes (T^{\mathfrak{s}|_2}_{y,j}T^{\mathfrak{b}|_2}_{x,l}-T^{\mathfrak{s}|_2}_{x,j}T^{\mathfrak{b}|_2}_{y,l})\\
&=\Delta(T^{\mathfrak{p}}_{k,j}T^{\mathfrak{q}}_{i,l}-T^{\mathfrak{q}}_{k,j}T^{\mathfrak{p}}_{i,l}),
\end{align*}
where the second equality is due to \eqref{gatt1}. We have proved the relation \eqref{gatt1}.

As for the relation \eqref{gatt2}, by a direct computation, we obtain
\begin{align*}
[\Delta(T^{(0)}_{i,j}),\Delta(T^{\mathfrak{q},\tau}_{k,l})]
&=([T^{(0)}_{i,j},T^{\mathfrak{q}|_1}_{k,l}])\otimes1+1\otimes [T^{(0)}_{i,j},T^{\mathfrak{q}|_2}_{k,l}]\\
&\quad+\hbar\sum_{y=1}^n\limits\sum_{\substack{\mathfrak{r}+\mathfrak{s}=\mathfrak{q},\\\mathfrak{r},\mathfrak{s}\neq(0)}}[T^{(0)}_{i,j},T^{\mathfrak{r}|_1}_{k,y}]\otimes T^{\mathfrak{s}|_2}_{y,l}+\hbar\sum_{y=1}^n\limits\sum_{\substack{\mathfrak{r}+\mathfrak{s}=\mathfrak{q},\\\mathfrak{r},\mathfrak{s}\neq(0)}}T^{\mathfrak{r}|_1}_{k,y}\otimes [T^{(0)}_{i,j},T^{\mathfrak{s}|_2}_{y,l}]\\
&=(\delta_{j,k}T^{\mathfrak{q}|_1}_{i,l}-\delta_{i,l}T^{\mathfrak{q}|_1}_{k,j})\otimes1+1\otimes(\delta_{j,k}T^{\mathfrak{q}|_2}_{i,l}-\delta_{i,l}T^{\mathfrak{q}|_2}_{k,j})\\
&\quad+\hbar\sum_{y=1}^n\limits\sum_{\substack{\mathfrak{r}+\mathfrak{s}=\mathfrak{q},\\\mathfrak{r},\mathfrak{s}\neq(0)}}(\delta_{j,k}T^{\mathfrak{r}|_1}_{i,y}-\delta_{y,i}T^{\mathfrak{r}|_1}_{k,j})T^{\mathfrak{s}|_2}_{y,l}+\hbar\sum_{y=1}^n\limits\sum_{\substack{\mathfrak{r}+\mathfrak{s}=\mathfrak{q},\\\mathfrak{r},\mathfrak{s}\neq(0)}}T^{\mathfrak{r}|_1}_{k,y}\otimes (\delta_{j,y}T^{\mathfrak{s}|_2}_{i,l}-\delta_{i,l}T^{\mathfrak{s}|_2}_{y,j})\\
&=\Delta(\delta_{j,k}T^{\mathfrak{q},\tau}_{i,l}-\delta_{i,l}T^{\mathfrak{q},\tau}_{k,j}),
\end{align*}
where the second equality is due to \eqref{gatt2}. We have proved the compatibility with \eqref{gatt2}.
\end{proof}
The following corollary can be proved by the same proof as Theorem~\ref{Thm1}.
\begin{Corollary}
There exists a homomorphism 
\[\Delta^{\text{op}}\colon Y_\hbar(n)\to Y_\hbar(n)\otimes Y_\hbar(n)\]
given by
\begin{align*}
\Delta^{\text{op}}(T^{\mathfrak{p}}_{i,j})&=T^{\mathfrak{p}|_2}_{i,j}\otimes1+1\otimes T^{\mathfrak{p}|_1}_{i,j}+\hbar\sum_{\substack{\mathfrak{r}+\mathfrak{s}=\mathfrak{p},\\\mathfrak{r},\mathfrak{s}\neq(0)}}T^{\mathfrak{r}|_2}_{i,x}\otimes T^{\mathfrak{s}|_1}_{x,j}.
\end{align*}
\end{Corollary}
Considering the result in Sections 4 and 6, it is expected that the Yangian associated with $X^*$ contains $\Aut(X^*)$. We extend the associative algebra $Y_\hbar(n)$
\begin{Lemma}
We can give an action of $\Aut(X^*)^\ast$ and $\Aut(X^*)$ on $Y_\hbar(n)$ by $\sigma T^{(x_w)}_{i,j}=T^{(x_{\sigma_0\sigma(w)})}_{i,j}$ and $\sigma T^{(x_w)}_{i,j}=T^{(x_{\sigma(w)})}_{i,j}$ respectively.
\end{Lemma}
\begin{proof}
We only show the case for $\Aut(X^*)^\ast$. The another case is trivial. For $\sigma,\tau\in\Aut(X^*),\mathfrak{p}\in P$ and $1\leq i,j\leq n$, we obtain
\begin{align*}
\sigma(\tau T^{\mathfrak{p}}_{i,j}))=T^{\sigma_0\sigma\tau_0\tau(\mathfrak{p})}_{i,j}=T^{\sigma_0\tau_0(\tau_0\sigma\tau_0)\tau(\mathfrak{p})}_{i,j}=T^{(\sigma\ast\tau)(\mathfrak{p}}_{i,j}=(\sigma\ast\tau)(T^{\mathfrak{p}}_{i,j}).
\end{align*}
\end{proof}
Let $Y^\ast(n)$ and $\widetilde{Y}(n)$ be $Y_\hbar(n)\rtimes\Aut(X^*)^*$ and $Y_\hbar(n)\rtimes\Aut(X^*)^\ast$.
\begin{Theorem}
There exists a homomorphism 
\[\widetilde{\Delta}\colon Y^\ast(n)\to \widetilde{Y}(n)\otimes \widetilde{Y}(n)\]
given by
\begin{align*}
\widetilde{\Delta}(\sigma)&=\sigma|_1\otimes\sigma|_2\otimes\sigma_0,\\
\widetilde{\Delta}(T^{\mathfrak{p}}_{i,j})&=T^{\mathfrak{p}|_1}_{i,j}\otimes1+1\otimes T^{\mathfrak{p}|_2}_{i,j}+\hbar\sum_{\substack{\mathfrak{r}+\mathfrak{s}=\mathfrak{p},\\\mathfrak{r},\mathfrak{s}\neq(0)}}T^{\mathfrak{r}|_1}_{i,x}\otimes T^{\mathfrak{s}|_2}_{x,j}.
\end{align*}
\end{Theorem}
\begin{proof}
By Lemma~\ref{Lem980} and Theorem~\ref{Thm1}, it is enough to show the compatibility with the relation
\[
\widetilde{\Delta}(\sigma)\widetilde{\Delta}(T^{\mathfrak{p}}_{i,j})=\widetilde{\Delta}(T^{\sigma_0\sigma(\mathfrak{p})}_{i,j})
\]
for $\sigma\in\Aut(X^*),\mathfrak{p}\in P$ and $1\leq i,j\leq n$. By a direct computation, we obtain
\begin{align*}
\widetilde{\Delta}(\sigma)\widetilde{\Delta}(T^{\mathfrak{p}}_{i,j})
&=\sigma|_1T^{\mathfrak{p}|_1}_{i,j}\otimes\sigma_21\otimes\sigma_01+\sigma|_11\otimes \sigma|_2T^{\mathfrak{p}|_2}_{i,j}\otimes\sigma_01+\hbar\sum_{\mathfrak{r}+\mathfrak{s}=\mathfrak{p},\mathfrak{r},\mathfrak{s}\neq(0)}\sigma|_1T^{\mathfrak{p}|_1}_{i,x}\otimes T^{\mathfrak{s}|_2}_{x,j}\otimes\sigma_01\\
&=T^{\sigma(\mathfrak{p})|_1}_{i,j}\otimes1\otimes1+1\otimes T^{\sigma(\mathfrak{p})|_2}_{i,j}\otimes1+\hbar\sum_{\mathfrak{r}+\mathfrak{s}=\mathfrak{p},\mathfrak{r},\mathfrak{s}\neq(0)}T^{\sigma(\mathfrak{p})|_1}_{i,x}\otimes T^{\sigma(\mathfrak{s})|_2}_{x,j}\otimes1\\
&=\Delta(T^{\sigma(\mathfrak{p})}_{i,j}).
\end{align*}
Thus, $\widetilde{\Delta}$ is well-defined.
\end{proof}
Let us connect the target of $\widetilde{\Delta}$ with $\text{Aut}(X^*)^*$.
\begin{Lemma}\label{Lem981}
We can define an action of $\mathbb{C}[S_2]\ltimes\Aut(X^*)^\ast$ on $Y(n)$ by
\[
(\omega,\sigma)T^{\mathfrak{p}}_{i,j}=T^{\omega\sigma_0\sigma(\mathfrak{p})}_{i,j}
\]
for $\omega\in S_2,\sigma\in\Aut(X^*),\mathfrak{p}\in P$ and $1\leq i,j\leq n$.
\end{Lemma}
\begin{proof}
For $\sigma,\tau\in\Aut(X^*),\alpha,\beta\in S_2,\mathfrak{p}\in P$ and $1\leq i,j\leq n$, by a direct computation, we obtain
\begin{align*}
(\alpha,\sigma)(\beta,\tau)=(\alpha\beta,(\beta\sigma\beta)\ast\tau)=(\alpha\beta,\tau_0\beta\sigma\beta\tau_0\tau).
\end{align*}
By a direct computation, we obtain
\begin{align*}
((\alpha,\sigma)(\beta,\tau))T^{\mathfrak{p}}_{i,j}=T^{\alpha\beta\sigma_0\tau_0\tau_0\beta\sigma\beta\tau_0\tau}_{i,j}=T^{\alpha\sigma_0\sigma\beta\tau_0\tau}_{i,j}
\end{align*}
for $\mathfrak{p}\in P$ and $1\leq i,j\leq n$. On the other hand, we have
\begin{align*}
(\alpha,\sigma)((\beta,\tau)T^{\mathfrak{p}}_{i,j})=(\alpha,\sigma)T^{\alpha\sigma_0\sigma\beta\tau_0\tau}_{i,j}=T^{\beta\tau_0\tau}_{i,j}.
\end{align*}
Thus, we obtain $(\alpha,\sigma)((\beta,\tau)T^{\mathfrak{p}}_{i,j})=((\alpha,\sigma)(\beta,\tau))T^{\mathfrak{p}}_{i,j}$.
\end{proof}
Let $\widehat{Y}(n)$ be $Y(n)\rtimes(\mathbb{C}[S_2]\ltimes\Aut(X^*)^\ast)$. By Lemam~\ref{Lem980} (3) and Lemma~\ref{Lem981}, we obtain
\[
\widehat{\Delta}\colon \widetilde{Y}(n)\to\widehat{Y}(n)\otimes\widehat{Y}(n).
\]
We expect that this homomorphism can be extended to $\widetilde{Y}(n)$ and corresponds to the wreath coproduct for the Yangian associated with $X^*$.

Actually, by Lemma~\ref{Lem3.3}, we can consider another homomorphism for $\widetilde{Y}(n)$.
\begin{Theorem}\label{DD}
We can define a homomorphism 
\[
\Delta^D\colon \widetilde{Y}(n)\to\widetilde{Y}(n)^{\otimes 4}\rtimes\mathbb{C}[S_2]
\]
given by
\begin{align*}
\Delta^D(\sigma)&=(\sigma|_{\sigma(1)}\otimes\sigma_{\sigma(2)}\otimes \sigma|_{\sigma(2)}\otimes\sigma_{\sigma(1)},\sigma_0),\\
\Delta^D(T^{\mathfrak{p}}_{i,j})&=(\Delta(T^{\mathfrak{p}}_{i,j})\otimes1\otimes1,1)+(1\otimes1\otimes\Delta^{\text{op}}(T^{\mathfrak{p}}_{i,j}),1)
\end{align*}
where we set an action of $\mathbb{C}[S_2]$ on $\widetilde{Y}(n)^{\otimes 4}$ by
\[
s_{1,2}(v_1\otimes v_2\otimes v_3\otimes v_4)=v_3\otimes v_4\otimes v_1\otimes v_2.
\]
\end{Theorem}
\begin{proof}
It is enough to show $\Delta^D(\sigma)\Delta^D(T^{\mathfrak{p}}_{i,j})=\Delta^D(T^{\sigma(\mathfrak{p})}_{i,j})$ for $1\leq i,j\leq n,\sigma\in\Aut(X^*)$ and $\mathfrak{p}\in P$. By a direct computation, we have
\begin{align*}
\Delta^D(\sigma)\Delta^D(T^{\mathfrak{p}}_{i,j})
&=(\sigma|_{\sigma(1)}\otimes\sigma_{\sigma(2)}\otimes \sigma|_{\sigma(2)}\otimes\sigma_{\sigma(1)},\sigma_0)(\Delta(T^{\mathfrak{p}}_{i,j})\otimes1\otimes1,1)+(1\otimes1\otimes\Delta^{\text{op}}(T^{\mathfrak{p}}_{i,j}),1)\\
&=(T^{\sigma|_1\mathfrak{p}|_1}_{i,j}\otimes1\otimes1,1)+(1\otimes T^{\sigma|_2\mathfrak{p}|_2}_{i,j}\otimes1\otimes1,1)\\
&\quad+\hbar\sum_{\substack{\mathfrak{r}+\mathfrak{s}=\mathfrak{p},\\\mathfrak{r},\mathfrak{s}\neq(0)}}(T^{\sigma|_1\mathfrak{r}|_1}_{i,x}\otimes T^{\sigma|_2\mathfrak{s}|_2}_{x,j}\otimes1\otimes1,1)\\
&\quad+(1\otimes1\otimes T^{\sigma|_2\mathfrak{p}|_2}_{i,j},1)+(1\otimes1\otimes1\otimes T^{\sigma|_1\mathfrak{p}|_1}_{i,j},1)\\
&\quad+\hbar\sum_{\substack{\mathfrak{r}+\mathfrak{s}=\mathfrak{p},\\\mathfrak{r},\mathfrak{s}\neq(0)}}(1\otimes1\otimes T^{\sigma|_2\mathfrak{r}|_2}_{i,x}\otimes T^{\sigma|_1\mathfrak{s}|_1}_{x,j},1)
\end{align*}
if $\sigma_0=1$ and
\begin{align*}
\Delta^D(\sigma)\Delta^D(T^{\mathfrak{p}}_{i,j})
&=(\sigma|_{\sigma(1)}\otimes\sigma_{\sigma(2)}\otimes \sigma|_{\sigma(2)}\otimes\sigma_{\sigma(1)},\sigma_0)(\Delta(T^{\mathfrak{p}}_{i,j})\otimes1\otimes1,1)+(1\otimes1\otimes\Delta^{\text{op}}(T^{\mathfrak{p}}_{i,j}),1)\\
&=(T^{\sigma|_2\mathfrak{p}|_2}_{i,j}\otimes1\otimes1,\sigma_0)+(1\otimes T^{\sigma|_1\mathfrak{p}|_1}_{i,j}\otimes1\otimes1,\sigma_0)\\
&\quad+\hbar\sum_{\substack{\mathfrak{r}+\mathfrak{s}=\mathfrak{p},\\\mathfrak{r},\mathfrak{s}\neq(0)}}(T^{\sigma|_2\mathfrak{r}|_2}_{i,x}\otimes T^{\sigma|_1\mathfrak{s}|_1}_{x,j}\otimes1\otimes1,\sigma_0)\\
&\quad+(1\otimes1\otimes T^{\sigma|_1\mathfrak{p}|_1}_{i,j},\sigma_0)+(1\otimes1\otimes1\otimes T^{\sigma|_2\mathfrak{p}|_2}_{i,j},\sigma_0)\\
&\quad+\hbar\sum_{\substack{\mathfrak{r}+\mathfrak{s}=\mathfrak{p},\\\mathfrak{r},\mathfrak{s}\neq(0)}}(1\otimes1\otimes T^{\sigma|_1\mathfrak{r}|_21}_{i,x}\otimes T^{\sigma|_2\mathfrak{s}|_2}_{x,j},\sigma_0)
\end{align*}
if $\sigma_0=s_{1,2}$. Thus, we obtain $\Delta^D(\sigma)\Delta^D(T^{\mathfrak{p}}_{i,j})=\Delta^D(T^{\sigma(\mathfrak{p})}_{i,j})$.
\end{proof}
Based on Theorem~\ref{DD}, we can consider $\widetilde{Y}(n)$ as one of candidates for the Yangian associated with $X^*$.

\end{document}